\documentclass[aos,seceqn,citesort,dvips]{arximspdf}
\doi{10.1214/09-AOS759}
\volume{38}
\issue{3}
\pubyear{2010}
\firstpage{1568}
\lastpage{1592}

\makeatletter

\newcommand{\dd}{{d}}

\newtheorem{theorem}{Theorem}[section]
\newtheorem{lem}{Lemma}[section]
\newtheorem{pro}{Proposition}[section]
\newproclaim{exmp}{Example}[section]
\newproclaim{rem}{Remark}[section]
\newproclaim{fact}{Fact}[section]

\makeatother

\begin{document}
\begin{frontmatter}

\title{Statistical analysis of $k$-nearest neighbor
collaborative recommendation}
\runtitle{$k$-nearest neighbor collaborative recommendation}

\begin{aug}
\author[A]{\fnms{G\'erard} \snm{Biau}\corref{}\ead[label=e1]{gerard.biau@upmc.fr}\ead[label=u1,url]{http://www.lsta.upmc.fr/biau.html}},
\author[B]{\fnms{Beno\^it} \snm{Cadre}\ead[label=e2]{Benoit.cadre@bretagne.ens-cachan.fr}\ead[label=u2,url]{http://w3.bretagne.ens-cachan.fr/math/people/benoit.cadre}} and
\author[C]{\fnms{Laurent} \snm{Rouvi\`ere}\ead[label=e3]{laurent.rouviere@ensai.fr}\ead[label=u3,url]{http://www.ensai.com/laurent-rouviere-rub,78.html}}
\runauthor{G. Biau, B. Cadre and L. Rouvi\`ere}
\affiliation{Universit\'e Paris VI, ENS Cachan-Bretagne and CREST-ENSAI}
\address[A]{G. Biau\\
LSTA and LPMA \\
Universit\'e Paris VI \\
Bo\^{\i}te 158, 175 rue du Chevaleret\\
75013 Paris\\
France\\
\printead{e1}\\
\printead{u1}}
\address[B]{B. Cadre\\
IRMAR, ENS Cachan Bretagne, CNRS, UEB \\
Campus de Ker Lann\\
Avenue Robert Schuman\\
35170 Bruz\\
France\\
\printead{e2}\\
\printead{u2}}
\address[C]{L. Rouvi\`ere \\
CREST-ENSAI, IRMAR, UEB\\
Campus de Ker Lann\\
Rue Blaise Pascal, BP 37203\\
35172 Bruz Cedex\\
France\\
\printead{e3}\\
\printead{u3}}
\end{aug}

\pdfauthor{Gerard Biau, Benoit Cadre, Laurent Rouviere}

% HISTORY:
\received{\smonth{3} \syear{2009}}
\revised{\smonth{9} \syear{2009}}

% ABSTRACT
%
\begin{abstract}
Collaborative recommendation is an information-filtering technique that
attempts to present information items that are likely of interest to an
Internet user. Traditionally, collaborative systems deal with
situations with two types of variables, users and items. In its most
common form, the problem is framed as trying to estimate ratings for
items that have not yet been consumed by a user. Despite wide-ranging
literature, little is known about the statistical properties of
recommendation systems. In fact, no clear probabilistic model even
exists which would allow us to precisely describe the mathematical
forces driving collaborative filtering. To provide an initial
contribution to this, we propose to set out a general sequential
stochastic model for collaborative recommendation. We offer an in-depth
analysis of the so-called cosine-type nearest neighbor collaborative
method, which is one of the most widely used algorithms in
collaborative filtering, and analyze its asymptotic performance as the
number of users grows. We establish consistency of the procedure under
mild assumptions on the model. Rates of convergence and examples are
also provided.
\end{abstract}

% KEYWORDS
%

\end{frontmatter}

%s1 ###
\section{Introduction}\label{sec1}

Collaborative recommendation is a Web information-fil\-tering technique
that typically gathers information about your personal interests and
compares your profile to other users with similar tastes. The goal of
this system is to give personalized recommendations, whether this be
movies you might enjoy, books you should read or the next restaurant
you should go to.

There has been much work done in this area over the past decade since
the appearance of the first papers on the subject in the
mid-90s \cite{HSRF95,RISBR94,SM95}. Stimulated by an abundance of practical
applications, most of the research activity to date has focused on
elaborating various heuristics and practical methods
\cite{BHK1998,HCMRK2000,SMH2007} so as to provide personalized
recommendations and help Web users deal with information overload.
Examples of such applications include recommending books, people,
restaurants, movies, CDs and news. Websites such as amazon.com,
match.com, movielens.org and allmusic.com already have recommendation
systems in operation. We refer the reader to the surveys by
\cite{AT2005} and \cite{ASST2005} for a broader picture of the field,
an overview of results and many related references.

Traditionally, collaborative systems deal with situations with two
types of variables, \textit{users} and \textit{items}. In its most common
form, the problem is framed as trying to estimate \textit{ratings} for
items that have \textit{not} yet been consumed by a user. The
recommendation process typically starts by asking users a series of
questions about items they liked or did not like. For example, in a
movie recommendation system, users initially rate some subset of films
they have already seen. Personal ratings are then collected in a matrix
where each row represents a user,
each column an item, and entries in the matrix represent a given user's
rating of a given item. An example is presented in Table \ref{tableau1}
where ratings are specified on a
scale from 1 to~10, and ``NA'' means that the user has not rated the
corresponding film.

%t1 ###
%
\begin{table}%[b]
\caption{A (subset of a) ratings matrix for a movie recommendation system.
Ratings are specified on a scale from $1$ to $10$, and ``NA'' means
that the user has not rated the corresponding film}
\label{tableau1}
\begin{tabular*}{\tablewidth}{@{\extracolsep{\fill}}lcccccc@{}}
\hline
& \textbf{Armageddon} & \textbf{Platoon} & \textbf{Rambo} & \textbf{Rio Bravo} &
\textbf{Star wars} & \textbf{Titanic} \\
\hline
{Jim} & NA & 6 & 7 & 8 & 9 & NA\\
{James} & 3 & NA & 10 & NA & 5 & 7\\
{Steve} & 7 & NA & 1 & NA & 6 & NA\\
{Mary} & NA & 7 & 1 & NA & 5 & 6\\
{John} & NA & 7 & NA & NA & 3 & 1\\
{Lucy} & 3 & 10 & 2 & 7 & NA & 4\\
{Stan} & NA & 7 & NA & NA & 1 & NA\\
{Johanna} & 4 & 5 & NA & 8 & 3 & 9\\
[4pt]
{Bob} & NA & 3 & 3 & 4 & 5 & $\mathbf{?}$\\
\hline
\end{tabular*}\vspace*{3pt}
\end{table}

Based on this prior information, the recommendation engine must be able
to automatically furnish ratings of as-yet unrated items and then
suggest appropriate recommendations based on these predictions. To do
this, a number of practical methods have been proposed, including
machine learning-oriented techniques \cite{ABEV2008}, statistical
approaches \cite{SKKR2001} and numerous other ad hoc rules
\cite{ASST2005}. The collaborative filtering issue may be viewed as a
special instance of the problem of inferring the many missing entries
of a data matrix. This field, which has very recently emerged, is known
as the matrix completion problem and
comes up in many areas of science and engineering, including
collaborative filtering, machine learning, control, remote sensing and
computer vision. We will not pursue this promising approach, and refer
the reader to \cite{candes1}
and \cite{candes2} who
survey the literature on matrix completion. These authors show in
particular that under suitable conditions, one can recover an unknown
low rank matrix from a nearly minimal set of entries by solving a
simple convex optimization problem.

In most of the approaches, the crux is to identify users whose
tastes/rat\-ings are ``similar'' to the user we would like to advise.
The similarity measure assessing proximity between users may vary
depending on the type of application
but is typically based on a correlation or cosine-type approach
\cite{SKKR2001}.

Despite wide-ranging literature, very little is known about the
statistical properties of recommendation systems. In fact, no clear
probabilistic model even exists allowing us to precisely describe the
mathematical forces driving collaborative filtering. To provide an
initial contribution to this, we propose in the present paper to set
out a general stochastic model for collaborative recommendation and
analyze its asymptotic performance as the number of users grows.

The document is organized as follows. In Section \ref{sec2}, we provide a
sequential stochastic model for collaborative recommendation and
describe the statistical problem. In the model we analyze, unrated
items are estimated by averaging ratings of users who are ``similar''
to the user we would like to advise. The similarity is assessed by a
cosine-type measure, and unrated items are estimated using a
\mbox{$k_n$-nearest} neighbor-type regression estimate which is indeed one of
the most widely used procedures in collaborative filtering. It turns
out that the choice of the cosine proximity as a similarity measure
imposes constraints on the model which are discussed in Section \ref{sec3}.
Under mild assumptions, consistency of the estimation procedure is
established in Section \ref{sec4} whereas rates of convergence are discussed in
Section \ref{sec5}. Illustrative examples are given throughout the document, and
proofs of some technical results are postponed to Section \ref{sec6}.

%s2 ###
\section{A model for collaborative recommendation}\label{sec2}

%s2.1 ###
\subsection{Ratings matrix and new users}\label{sec21}

Suppose that there are $d+1$ ($d \geq1$) possible items, $n$ users in
the ratings matrix (i.e., the database) and that users' ratings take
values in the set $(\{0\} \cup[1,s])^{d+1}$. Here, $s$ is a real
number greater than 1 corresponding to the maximal rating, and, by
convention, the symbol 0 means that the user has not rated the item
(same as ``NA''). Thus the ratings matrix has $n$ rows, $d+1$ columns
and entries from $\{0\} \cup[1,s]$. For example, $n=8$, $d=5$ and
$s=10$ in Table \ref{tableau1} which will be our toy example throughout
this section. Then a new user, Bob, reveals some of his preferences for
the first time, rating some of the first $d$ items but \textit{not} the
$(d+1)$th (the movie \textit{Titanic} in Table \ref{tableau1}). We want
to design a strategy to predict Bob's rating of \textit{Titanic} using:
(i) Bob's ratings of some (or all) of the other $d$ movies and (ii)
the ratings matrix. This is illustrated in Table \ref{tableau1}, where
Bob has rated 4 out of the 5 movies.

The first step in our approach is to model the preferences of the new
user, Bob, by a random vector $(\mathbf{X},Y)$ of size $d+1$ taking
values in
the set $[1,s]^d \times[1,s]$. Within this framework, the random
variable $\mathbf{X}=(X_1,\ldots, X_d)$ represents Bob's preferences
pertaining to the first $d$ movies whereas $Y$, the (unobserved)
variable of interest, refers to the movie \textit{Titanic}. In fact, as
Bob does not necessarily reveals all his preferences at once, we do not
observe the variable $\mathbf{X}$, but instead some ``masked'' version
of it
denoted hereafter by $\mathbf{X}^{\star}$. The random variable
$\mathbf{X}^{\star}=(X^{\star}_1,\ldots, X^{\star}_d)$ is
naturally defined
by
\[
X_j^{\star}=
\cases{
X_j, &\quad if $j \in M$,\cr
0, &\quad otherwise,}
\]
where $M$ stands for some nonempty random subset of $\{1,\ldots, d\}$
indexing the movies which have been rated by Bob. Observe that the
random variable $\mathbf{X}^{\star}$ takes values in $(\{0\}\cup[1,s])^{d}$
and that $\|\mathbf{X}^{\star}\| \geq1$ where $\|\cdot\|$ denotes the usual
Euclidean norm on $\mathbb R^d$. In the example of Table \ref
{tableau1}, $M=\{2,3,4,5\}$ and (the realization of) $\mathbf
{X}^{\star}$ is
$(0,3,3,4,5)$.

We follow the same approach to model preferences of users already in
the database (Jim, James, Steve, Mary, etc. in Table \ref{tableau1}),
who will therefore be represented by a sequence of independent $[1,s]^d
\times[1,s]$-valued random pairs $(\mathbf{X}_1,Y_1),\ldots,
(\mathbf{X}_n,Y_n)$
from the distribution $(\mathbf{X},Y)$. A first idea for dealing with
potential nonresponses of a user $i$ in the ratings matrix ($i=1,\ldots
, n$) is to consider in place of $\mathbf{X}_i=(X_{i1},\ldots,
X_{id})$, its
masked version $\widetilde{\mathbf{X}}_i=(\widetilde{X}_{i1},\ldots,
\widetilde
{X}_{id})$ defined by
%
%e2.1 ###
%
\begin{equation}
\label{crise}
\widetilde{X}_{ij}=
\cases{
X_{ij}, &\quad if $j \in M_i \cap M$,\cr
0, &\quad otherwise,}
\end{equation}
where each $M_i$ is the random subset of $\{1,\ldots,d\}$ indexing the
movies which have been rated by user $i$. In other words, we only keep
in $\mathbf{X}_i$ items corated by both user $i$ \textit{and} the new
user---items which have not been rated by $\mathbf{X}$ and $\mathbf{X}_i$
are declared
noninformative and simply thrown away.

However, this model, which is static in nature, does not allow to take
into account the fact that, as time goes by, each user in the database
may reveal more and more preferences. This will, for instance,
typically be the case in the movie recommendation system of Table \ref
{tableau1} where regular customers will update their ratings each time
they have seen a new movie. Consequently, model (\ref{crise}) is not
fully satisfying and must therefore be slightly modified to better
capture the sequential evolution of ratings.

%s2.2 ###
\subsection{A sequential model}\label{sec22}

A possible dynamical approach for collaborative recommendation is based
on the following protocol: users enter the database one after the other
and update their list of ratings sequentially in time. More precisely,
we suppose that at each time $i=1, 2,\ldots,$ a new user enters the
process and reveals his preferences for the first time while the $i-1$
previous users are allowed to rate new items. Thus, at time 1, there is
only one user in the database (Jim in Table \ref{tableau1}), and the
(nonempty) subset of items he decides to rate is modeled by a random
variable $M_1^1$ taking values in $\mathcal P^{\star}(\{1,\ldots, d\}
)$, the set of nonempty subsets of $\{1,\ldots,d\}$. At time 2, a new
user (James) enters the game and reveals his preferences according to a
$\mathcal P^{\star}(\{1,\ldots, d\})$-valued random variable $M_2^1$,
with the same distribution as $M_1^1$. At the same time, Jim (user 1)
may update his list of preferences, modeled by a random variable
$M_1^2$ satisfying $M_1^1\subset M_1^2$. The latter requirement just
means that the user is allowed to rate new items but not to remove his
past ratings. At time~3, a new user (Steve) rates items according to a
random variable $M_3^1$ distributed as~$M_1^1$, while user 2 updates
his preferences according to $M_2^2$ (distributed as $M_1^2$) and
user~1 updates his own according to $M_1^3$, and so on. This sequential
mechanism is summarized in Table \ref{tableau2}.

%t2 ###
%
\begin{table}
\tablewidth=250pt
\caption{A sequential model for preference updating}
\label{tableau2}
\begin{tabular*}{\tablewidth}{@{\extracolsep{\fill}}lcccccc@{}}
\hline
& \textbf{Time 1} & \textbf{Time 2} & $\bolds\cdots$ & \textbf{Time} $\bolds i$ & $\bolds\cdots$ & \textbf{Time} $\bolds n$\\
\hline
User 1 & $M_1^1$ & $M_1^2$ & $\cdots$ & $M_1^i $ & $\cdots$ & $M_1^n$
\\[3pt]
User 2 & & $M_2^1$ & $\cdots$ & $M_2^{i-1} $ & $\cdots$ & $M_2^{n-1}$
\\
\vdots& &&$\ddots$&$\vdots$& $\vdots$& $\vdots$\\
User $i$ & && & $M_i^1$ & $\cdots$ & $M_i^{n+1-i}$\\
\vdots& &&&& $\ddots$ &$\vdots$\\
User $n$ & &&&&&$M_n^1$\\
\hline
\end{tabular*}
\end{table}

By repeating this procedure, we end up at time $n$ with an upper
triangular array $(M_i^j)_{1\leq i \leq n, 1\leq j \leq n+1-i}$ of
random variables. A row in this array consists of a collection
$M_{i}^j$ of random variables for a given value of $i$, taking values
in $\mathcal P^{\star}(\{1,\ldots, d\})$ and satisfying the constraint
$M_{i}^j \subset M_{i}^{j+1}$. For a fixed $i$, the sequence $M_{i}^1
\subset M_{i}^{2}\subset\cdots$ describes the (random) way user $i$
sequentially reveals his preferences over time. Observe that the later
inclusions are not necessarily strict, so that a single user is not
forced to rate one more item at every single step.

Throughout the paper, we will assume that, for each $i$, the
distribution of the sequence of random variables $(M_i^n)_{n\geq1}$ is
independent of $i$, and is therefore distributed as a generic random
sequence denoted $(M^n)_{n\geq1}$, satisfying $M^1\neq\varnothing$ and
$M^n \subset M^{n+1}$ for all $n \geq1$. For the sake of coherence, we
assume that $M^1$ and $M$ [see (\ref{crise})] have the same distribution;
that is, the new abstract user $\mathbf{X}^{\star}$ may be regarded
as a
user entering the database for the first time. We will also suppose
that there exists a positive random integer $n_0$ such that
$M^{n_0}=\{1,\ldots, d\}$, and, consequently, $M^n=\{1,\ldots, d\}$
for all $n \geq n_0$. This requirement means that each user rates all
$d$ items after a (random) period of time. Last, we will assume that
the pairs $(\mathbf{X}_i,Y_i)$, $i=1,\ldots, n$, the sequences $(M_1^n)_{n
\geq1}$, $(M_2^n)_{n \geq1}, \ldots$ and the random variable $M$ are
mutually independent. We note that this implies that the users' ratings
are independent.

With this sequential point of view, improving on (\ref{crise}), we let
the masked version $\mathbf{X}_i^{(n)}=(X_{i1}^{(n)},\ldots,
X_{id}^{(n)})$ of
$\mathbf{X}_i$ be defined as
\[
X_{ij}^{(n)}=
\cases{
X_{ij}, &\quad if $j \in M_{i}^{n+1-i} \cap M$,\cr
0, &\quad otherwise.}
\]
Again, it is worth pointing out that, in the definition of $\mathbf{X}
_i^{(n)}$, items which have not been corated by both $\mathbf{X}$ \textit{and}
$\mathbf{X}_i$ are deleted. This implies in particular that $\mathbf
{X}_i^{(n)}$ may
be equal to $\mathbf{0}$, the $d$-dimensional null vector (whereas $\|
\mathbf{X}
^{\star}\|\geq1$ by construction).

Finally, in order to deal with possible nonanswers of database users
regarding the variable of interest (\textit{Titanic} in our movie
example), we introduce $(\mathcal R_n)_{n \geq1}$, a sequence of
random variables taking values in $\mathcal P^{\star}(\{1,\ldots, n\}
)$, such that $\mathcal R_n$ is independent of $M$ and the sequences
$(M_i^{n})_{n \geq1}$, and satisfying $\mathcal R_n \subset\mathcal
R_{n+1}$ for all $n \geq1$. In this formalism, $\mathcal R_n$
represents the subset, which is assumed to be nonempty, of users who
have already provided information about \textit{Titanic} at time~$n$.
For example, in Table \ref{tableau1}, only James, Mary, John, Lucy and
Johanna have rated \textit{Titanic} and therefore (the realization of)
$\mathcal R_n$ is $\{2,4,5,6,8\}$.

%s2.3 ###
\subsection{The statistical problem}\label{sec23}

To summarize the model so far, we have at hand at time $n$ a sample of
random pairs $(\mathbf{X}_1^{(n)}, Y_1),\ldots, (\mathbf
{X}_n^{(n)},Y_n)$ and our
mission is to predict the score $Y$ of a new user represented by
$\mathbf{X}
^\star$. The variables $\mathbf{X}_1^{(n)},\ldots, \mathbf
{X}_n^{(n)}$ model the
database users' revealed preferences with respect to the first $d$
items. They take values in $(\{0\}\cup[1,s])^d$, where a 0 at
coordinate $j$ of $\mathbf{X}_i^{(n)}$ means that the $j$th product
has not
been corated by both user $i$ and the new user. The variable $\mathbf{X}
^\star
$ takes values in $(\{0\}\cup[1,s])^d$ and satisfies $\| \mathbf
{X}^{\star
}\|
\geq1$. The random variables $Y_1,\ldots, Y_n$ model users' ratings of
the product of interest. They take values in $[1,s]$ and, at time $n$,
we only see a nonempty (random) subset of $\{Y_1,\ldots, Y_n\}$,
indexed by $\mathcal R_n$.

The statistical problem with which we are faced is to estimate the
regression function $\eta(\mathbf{x}^{\star})=\mathbb E[Y|\mathbf{X}
^{\star}=\mathbf{x}
^{\star}]$. For this goal, we may use the database observations
$(\mathbf{X}
_1^{(n)}, Y_1),\ldots, (\mathbf{X}_n^{(n)},Y_n)$ in order to
construct an
estimate $\eta_n(\mathbf{x}^{\star})$ of $\eta(\mathbf{x}^{\star
})$. The approach we
explore in this paper is a cosine-based $k_n$-nearest neighbor
regression method, one of the most widely used algorithms in
collaborative filtering (see, e.g., \cite{SKKR2001}).

Given $\mathbf{x}^{\star} \in(\{0\}\cup[1,s])^d-\mathbf{0}$ and the
sample $(\mathbf{X}
_1^{(n)}, Y_1),\ldots, (\mathbf{X}_n^{(n)},Y_n)$, the idea of the cosine-type
$k_n$-nearest neighbor (NN) regression method is to estimate $\eta
(\mathbf{x}
^{\star})$ by a local averaging over those $Y_i$ for which: (i)
$\mathbf{X}
_i^{(n)}$ is ``close'' to $\mathbf{x}^\star$, and (ii) $i \in
\mathcal R_n$,
that is, we effectively ``see'' the rating $Y_i$. For this, we scan
through the $k_n$ neighbors of $\mathbf{x}^{\star}$ among the
database users
$\mathbf{X}_i^{(n)}$ for which $i \in\mathcal R_n$ and estimate $\eta
(\mathbf{x}
^{\star})$ by averaging the $k_n$ corresponding $Y_i$. The closeness
between users is assessed by a cosine-type similarity, defined for
$\mathbf{x}
=(x_1,\ldots, x_d)$ and $\mathbf{x}'=(x'_1,\ldots, x'_d)$ in $(\{0\}
\cup
[1,s])^d$ by
\[
\bar S(\mathbf{x}, \mathbf{x}')=\frac{\sum_{j\in\mathcal
J}x_jx_j'}{\sqrt{\sum_{j\in
\mathcal J}x_j^2}\sqrt{\sum_{j\in\mathcal J}x_j'^2}},
\]
where $\mathcal J=\{j\in\{1,\ldots,d\}\dvtx x_j\neq0$ and
$x_j'\neq
0\}$, and, by convention, $\bar S(\mathbf{x}, \mathbf{x}')=0$ if
$\mathcal J=\varnothing
$. To understand the rationale behind this proximity measure, just note
that if $\mathcal J=\{1,\ldots,d\}$ then $\bar S(\mathbf{x}, \mathbf
{x}')$ coincides
with $\cos(\mathbf{x}, \mathbf{x}')$; that is, two users are
``close'' with respect
to $\bar S$ if their ratings are more or less proportional. However,
the similarity $\bar S$, which will be used to measure the closeness
between $\mathbf{X}^\star$ (the new user) and $\mathbf{X}_i^{(n)}$
(a database user)
ignores possible nonanswers in $\mathbf{X}^{\star}$ or $\mathbf
{X}_i^{(n)}$, and is
therefore more adapted to the recommendation setting. For example, in
Table \ref{tableau1},
\begin{eqnarray*}
\bar S(\mbox{Bob}, \mbox{Jim})&=&\bar S ((0,3,3,4,5),(0,6,7,8,9))
\\
&=&\bar
S((3,3,4,5),(6,7,8,9))
 \approx0.99,
\end{eqnarray*}
whereas
\begin{eqnarray*}
\bar S(\mbox{Bob},\mbox{Lucy})&=&\bar S((0,3,3,4,5),(3,10,2,7,0)) \\
&=&\bar
S((3,3,4),(10,2,7))
\approx0.89.
\end{eqnarray*}
Next, fix $\mathbf{x}^{\star} \in(\{0\}\cup[1,s])^d-\mathbf{0}$, and suppose
for simplification that $M\subset M_{i}^{n+1-i}$ for each $i\in
\mathcal R_n$. In this case, it is easy to see that $\mathbf
{X}_i^{(n)}= \mathbf{X}
_i^{\star}=(X_{i1}^{\star},\ldots,X_{id}^{\star})$ where
\[
X_{ij}^{\star}=
\cases{
X_{ij}, &\quad if $j \in M$,\cr
0, &\quad otherwise.}
\]
Besides, $Y_i \geq1$,
%
%e2.2 ###
%
\begin{equation}
\label{isup}
\bar S(\mathbf{x}^{\star},\mathbf{X}_i^{\star})=\cos(\mathbf
{x}^{\star
},\mathbf{X}_i^{\star})>0
\end{equation}
and an elementary calculation shows that the positive real number $y$
which maximizes the similarity between $(\mathbf{x}^{\star},y)$ and
$(\mathbf{X}
_i^{\star},Y_i)$, that is,
\[
\bar S ((\mathbf{x}^{\star},y),(\mathbf{X}_i^{\star},Y_i) )=\frac
{\sum
_{j\in
M}x_j^\star X_{ij}^\star+yY_i}{\sqrt{\sum_{j\in M}{x_j^\star
}^2+y^2}\sqrt{\sum_{j\in M}{X_{ij}^\star}^2+Y_i^2}},
\]
is given by
\[
y=\frac{\|\mathbf{x}^{\star}\|}{\|\mathbf{X}_i^{\star}\| \cos
(\mathbf
{x}^{\star},\mathbf{X}_i^{\star
})} Y_i.
\]
This suggests the following regression estimate $\eta_n(\mathbf
{x}^{\star})$
of $\eta(\mathbf{x}^{\star})$:
%
%e2.3 ###
%
\begin{equation}
\label{benoitcadreboittrop}
\eta_n(\mathbf{x}^{\star})=\| \mathbf{x}^{\star} \| \sum_{i \in
\mathcal R_n}
W_{ni}(\mathbf{x}^{\star}) \frac{Y_i}{\|\mathbf{X}_i^{\star} \|},
\end{equation}
where the integer $k_n$ satisfies $1\leq k_n \leq n$ and
\[
W_{ni}(\mathbf{x}^{\star})=
\cases{
1/k_n, &\quad if $\mathbf{X}_i^{\star}$ is among the $k_n$-MS of $\mathbf{x}
^{\star}$ in $\{\mathbf{X}_i^{\star}, i \in\mathcal R_n\}$,\cr
0, &\quad otherwise.}
\]
In the above definition, the acronym ``MS'' (for most similar) means
that we are searching for the $k_n$ ``closest'' points of $\mathbf
{x}^{\star}$
within the set $\{\mathbf{X}_i^{\star}, i \in\mathcal R_n\}$ using the
similarity $\bar S$---or, equivalently here, using the cosine
proximity [by identity (\ref{isup})]. Note that the cosine term has
been removed since it has
asymptotically no influence on the estimate, as can be seen by a slight
adaptation of the arguments of the proof of Lemma 6.1, Chapter 6 in
\cite{gykokrha}. The estimate $\eta_n(\mathbf{x}^{\star})$ is
called the \textit{cosine-type $k_n$-NN regression estimate} in the collaborative
filtering literature. Now, recalling that definition
(\ref{benoitcadreboittrop}) makes sense only when $M\subset
M_{i}^{n+1-i}$ for each $i \in\mathcal R_n$ (that is, $\mathbf{X}_i^{(n)}=
\mathbf{X}_i^{\star}$), the next step is to extend the definition of
$\eta_n(\mathbf{x}^{\star})$ to the general case. In view of
(\ref{benoitcadreboittrop}), the most natural approach is to simply put
%
%e2.4 ###
%
\begin{equation}
\label{benoitcadreboitvraimenttrop}
\eta_n(\mathbf{x}^{\star})=\| \mathbf{x}^{\star} \| \sum_{i \in
\mathcal R_n}
W_{ni}(\mathbf{x}^{\star}) \frac{Y_i}{\|\mathbf{X}_i^{(n)} \|},
\end{equation}
where
\[
W_{ni}(\mathbf{x}^{\star})=
\cases{
1/k_n, &\quad if $\mathbf{X}_i^{(n)}$ is among the $k_n$-MS of $\mathbf{x}
^{\star}$ in $\bigl\{\mathbf{X}_i^{(n)}, i \in\mathcal R_n\bigr\}$,\cr
0, &\quad otherwise.}
\]
The acronym ``MS'' in the weight $W_{ni}(\mathbf{x}^{\star})$ means
that the
$k_n$ closest database points of $\mathbf{x}^{\star}$ are computed according
to the similarity
\[
S \bigl(\mathbf{x}^\star,\mathbf{X}_i^{(n)} \bigr)=p_i^{(n)}\bar S \bigl(\mathbf
{x}^\star
,\mathbf{X}
_i^{(n)} \bigr)\qquad \mbox{with } p_i^{(n)}=\frac{|{M_i^{n+1-i}\cap
M}|}{|{M}|},
\]
(here and throughout, notation $|{A}|$ means the cardinality of the
finite set $A$). The factor $p_i^{(n)}$ in front of $\bar S$ is a
penalty term which, roughly, avoids over promotion of the last users
entering the database. Indeed, the effective number of items rated by
these users will be eventually low, and, consequently, their $\bar
S$-proximity to $\mathbf{x}^{\star}$ will tend to remain high. On the other
hand, for fixed $i$ and $n$ large enough, we know that $M\subset
M_{i}^{n+1-i}$ and $\mathbf{X}_i^{(n)}=\mathbf{X}_i^{\star}$. This implies
$p_i^{(n)}=1$, $S(\mathbf{x}^\star,\mathbf{X}_i^{(n)})=\bar
S(\mathbf
{x}^\star,\mathbf{X}_i^{\star
})=\cos(\mathbf{x}^\star,\mathbf{X}_i^{\star})$ and shows that definition
(\ref
{benoitcadreboitvraimenttrop}) generalizes definition (\ref
{benoitcadreboittrop}). Therefore, we take the liberty to still call
the estimate (\ref{benoitcadreboitvraimenttrop}) the cosine-type
$k_n$-NN regression estimate.
\begin{rem} A smoothed version of the similarity $S$ could also be
considered, typically,
\[
S \bigl(\mathbf{x}^\star,\mathbf{X}_i^{(n)} \bigr)=\psi\bigl(p_i^{(n)} \bigr)\bar
S \bigl(\mathbf{x}^\star,\mathbf{X}_i^{(n)} \bigr),
\]
where $\psi\dvtx[0,1]\to[0,1]$ is a nondecreasing map satisfying
$\psi(1/2)<1$ (assuming $|{M}|\geq2$). For example, the choice
$\psi
(p)=\sqrt{p}$ tends to promote users with a low number of rated items,
provided the items corated by the new user are quite similar. In the
present paper, we shall only consider the case $\psi(p)=p$, but the
whole analysis carries over without difficulties for general functions
$\psi$.
\end{rem}
\begin{rem}
Another popular approach to measure the closeness between users is the
Pearson correlation coefficient. The extension of our results to
Pearson-type similarities is not straightforward and more work is
needed to address this challenging question. We refer the reader to
\cite{CKJ06} and \cite{MLR03} for a comparative study and comments on
the choice of the similarity.
\end{rem}

Finally, for definiteness of the estimate $\eta_n(\mathbf{x}^{\star
})$, some
final remarks are in order:
\begin{longlist}[(iii)]
\item[(i)] If\vspace*{-2pt} $S(\mathbf{x}^{\star},\mathbf{X}_i^{(n)})=S(\mathbf
{x}^{\star},\mathbf{X}
_j^{(n)})$, i.e., $\mathbf{X}_i^{(n)}$ and $\mathbf{X}_j^{(n)}$ are
equidistant from
$\mathbf{x}^{\star}$, then we have a tie, and, for example, $\mathbf{X}
_i^{(n)}$ may
be declared ``closer'' to $\mathbf{x}^\star$ if $i<j$; that is, tie-breaking
is done by indices.
\item[(ii)] If $|\mathcal R_n|<k_n$, then the weights
$W_{ni}(\mathbf{x}
^{\star})$ are not defined. In this case, we conveniently set
$W_{ni}(\mathbf{x}^{\star})=0$; that is, $\eta_n(\mathbf{x}^{\star})=0$.
\item[(iii)] If $\mathbf{X}_i^{(n)}=\mathbf0$, then we take
$W_{ni}(\mathbf{x}
^{\star})=0$, and we adopt the convention $0\times\infty=0$ for the
computation of $\eta_n(\mathbf{x}^{\star})$.
\item[(iv)] With the above conventions, the identity $\sum_{i \in
\mathcal R_n}W_{ni}(\mathbf{x}^{\star})\leq1$ holds in each case.
\end{longlist}

%
%s3 ###
\section{The regression function}\label{sec3}

Our objective in Section \ref{sec4} will be to establish consistency of the
estimate $\eta_n(\mathbf{x}^{\star})$ defined in (\ref
{benoitcadreboitvraimenttrop}) toward the regression function~$\eta
(\mathbf{x}^{\star})$. To reach this goal, we first need to analyze the
properties of $\eta(\mathbf{x}^{\star})$. Surprisingly, the special form
of $\eta_n(\mathbf{x}^{\star})$ constrains the shape of $\eta
(\mathbf{x}^{\star})$.
This is stated in Theorem \ref{forme} below.
\begin{theorem}
\label{forme}
Suppose that $\eta_n(\mathbf{X}^\star) \to\eta(\mathbf{X}^{\star})$ in
probability as $n \to\infty$. Then
\[
\eta(\mathbf{X}^\star)=\|\mathbf{X}^\star\| \mathbb E \biggl[ \frac{Y}{\|
\mathbf{X}^\star\|} \bigg| \frac{\mathbf{X}^\star}{\|\mathbf{X}^\star\|}
\biggr] \qquad\mbox{a.s.}
\]
\end{theorem}
\begin{pf}
Recall that
\[
\eta_n(\mathbf{X}^{\star})=\|\mathbf{X}^{\star}\|\sum_{i\in\mathcal R_n}
W_{ni}(\mathbf{X}^\star) \frac{Y_i}{\|\mathbf{X}_i^{(n)}\|}
\]
and let
\[
\varphi_n(\mathbf{X}^\star)=\sum_{i\in\mathcal R_n} W_{ni}(\mathbf{X}^\star)
\frac{Y_i}{\|\mathbf{X}_i^{(n)}\|}.
\]
Since $(\eta_n(\mathbf{X}^\star))_n$ is a Cauchy sequence in probability
and $\|\mathbf{X}^\star\|\geq1$, the sequence $(\varphi_n(\mathbf{X}^\star
))_n$ is also a Cauchy sequence. Thus there exists a measurable
function $\varphi$ on $\mathbb{R}^d$ such that $\varphi_n(\mathbf{X}^\star
)\to\varphi(\mathbf{X}^\star)$ in probability. Using the fact that $0
\leq\varphi_n(\mathbf{X}^\star) \leq s$ for all $n\geq1$, we conclude
that $0 \leq\varphi(\mathbf{X}^\star)\leq s$ a.s. as well.

Let us extract a sequence $(n_k)_k$ satisfying
$\varphi_{n_k}(\mathbf{X}^\star)\to\varphi(\mathbf{X}^\star)$ a.s.
Observing that, for $\mathbf{x}^\star\neq\mathbf0$,
\[
\varphi_{n_k}(\mathbf{x}^\star)=\varphi_{n_k} \biggl(\frac{\mathbf{x}^\star}{\|
\mathbf{x}^\star\|} \biggr),
\]
we may write
$\varphi(\mathbf{X}^\star)=\varphi(\mathbf{X}^\star/\|\mathbf{X}^\star\|)$
a.s. Consequently, the limit in probability of $(\eta
_n(\mathbf{X}^\star))_n$ is
\[
\|\mathbf{X}^\star\| \varphi\biggl(\frac{\mathbf{X}^\star}{\|\mathbf{X}^\star\|
} \biggr).
\]
Therefore, by the uniqueness of the limit, $\eta(\mathbf X^{\star})=\|
\mathbf{X}^\star\| \varphi(\mathbf{X}^\star/\|\mathbf{X}^\star\|)$ a.s. Moreover,
\begin{eqnarray*}
\varphi\biggl(\frac{\mathbf{X}^\star}{\|\mathbf{X}^\star\|} \biggr)
&=&
\mathbb E \biggl[\varphi\biggl(\frac{\mathbf{X}^\star}{\|\mathbf{X}^\star\|
} \biggr) \bigg|\frac{\mathbf{X}^\star}{\|\mathbf{X}^\star\|} \biggr]\\
&=& \mathbb E \biggl[\frac{\eta(\mathbf{X}^\star)}{\|\mathbf{X}^\star\|}
\bigg|\frac{\mathbf{X}^\star}{\|\mathbf{X}^\star\|} \biggr]\\
&=& \mathbb E \biggl[ \mathbb E \biggl[\frac{Y}{\|\mathbf{X}^\star\|}
\bigg| \mathbf{X}^\star\biggr] \bigg| \frac{\mathbf{X}^\star}{\|\mathbf{X}^\star\|}
\biggr]\\
&=& \mathbb E \biggl[ \frac{Y}{\|\mathbf{X}^\star\|} \bigg| \frac{\mathbf{X}^\star}{\|\mathbf{X}^\star\|} \biggr],
\end{eqnarray*}
since $\sigma(\mathbf{X}^\star/\|\mathbf{X}^\star\|)\subset\sigma(\mathbf{X}^\star)$. This completes the proof of the theorem.
\end{pf}

An important consequence of Theorem \ref{forme} is that if we intend to
prove any consistency result regarding the estimate $\eta_n(\mathbf
x^{\star})$, then we have to assume that the regression function $\eta
({\mathbf x}^{\star})$ has the special form
{\renewcommand{\theequation}{$\mathbf F$}
\begin{equation}
\label{eqF}
\eta(\mathbf{x}^\star)=\|\mathbf{x}^\star\|
\varphi(\mathbf{x}^\star)\qquad
\mbox{where } \varphi(\mathbf{x}^\star)= \mathbb E \biggl[ \frac{Y}{\|\mathbf{X}^\star
\|} \bigg| \frac{\mathbf{X}^\star}{\|\mathbf{X}^\star\|}=\frac{\mathbf{x}^\star
}{\|\mathbf{x}^\star\|} \biggr].
\end{equation}}
This will be our fundamental requirement throughout the paper, and it
will be denoted by (\ref{eqF}). In particular, if $\tilde\mathbf
{x}^\star
=\lambda\mathbf{x}^\star$ with $\lambda>0$, then $\eta(\tilde
\mathbf{x}^\star
)=\lambda\eta(\mathbf{x}^\star)$. That is, if two ratings $\mathbf
{x}^\star$ and
$\tilde\mathbf{x}^\star$ are proportional, then so must be the
values of the
regression function at $\mathbf{x}^\star$ and $\tilde\mathbf
{x}^\star$, respectively.

%s4 ###
\section{Consistency}\label{sec4}

In this section, we establish the $L_1$ consistency of the regression
estimate $\eta_n(\mathbf{x}^{\star})$ toward the regression function
$\eta(\mathbf{x}
^{\star})$. Using $L_1$ consistency is essentially a matter of taste,
and all the subsequent results may be easily adapted to $L_p$ norms
without too much effort. In the proofs, we will make repeated use of
the two following facts. Recall that, for a fixed $i \in\mathcal R_n$,
the random variable $\mathbf{X}_i^{\star}=(X_{i1}^{\star},\ldots
,X_{id}^{\star
})$ is defined by
\[
X_{ij}^{\star}=
\cases{
X_{ij}, &\quad if $j \in M$,\cr
0, &\quad otherwise,}
\]
and $\mathbf{X}_i^{(n)}= \mathbf{X}_i^{\star}$ as soon as $M\subset
M_{i}^{n+1-i}$.
Recall also that, by definition, $\|\mathbf{X}_i^{\star}\|\geq1$.
\begin{fact}
\label{FAIT1}
For each $i \in\mathcal R_n$,
\[
S(\mathbf{X}^{\star},\mathbf{X}_i^{\star})=\bar S(\mathbf
{X}^{\star},\mathbf{X}_i^{\star
})=\cos(\mathbf{X}
^{\star},\mathbf{X}_i^{\star})=1-\frac{1}{{2}} \dd^2 \biggl( \frac
{\mathbf{X}^{\star
}}{\|\mathbf{X}^{\star}\|},\frac{\mathbf{X}_i^{\star}}{\|\mathbf
{X}_i^{\star}\|} \biggr),
\]
where $\dd$ is the usual Euclidean distance on $\mathbb R^d$.
\end{fact}
\begin{fact}
\label{FAIT2}
Let, for all $i\geq1$,
\[
T_i=\min(k\geq i\dvtx M_i^{k+1-i}\supset M)
\]
be the first time instant when user $i$ has rated all the films indexed
by $M$. Set
%
%e4.1 ###
%
\begin{equation}
\label{eq:defLn}
\mathcal L_n=\{i\in\mathcal R_n\dvtx T_i\leq n\}
\end{equation}
and define, for $i\in\mathcal L_n$,
\[
W_{ni}^\star(\mathbf{x}^{\star})=
\cases{
1/k_n, &\quad if $\mathbf{X}_i^\star\mbox{ is among the $k_n$-MS of
}\mathbf{x}^{\star}$ in $\{\mathbf{X}_i^\star, i\in\mathcal L_n\}$, \cr
0, &\quad otherwise.}
\]
Then
\[
W_{ni}^\star(\mathbf{x^\star})=\cases{
1/k_n, &\quad if $\dfrac{\mathbf{X}_i^\star}{\|\mathbf{X}_i^\star\|}$ is among the $k_n$-NN of $\dfrac{\mathbf{x^{\star
}}}{\|{\mathbf{x}^{\star}}\|}$ in $ \biggl\{\dfrac{\mathbf{X}_i^\star}{\|
\mathbf{X}_i^\star\|},i\in\mathcal L_n \biggr\}$,\cr
0, &\quad otherwise,}
\]
where the $k_n$-NN are evaluated with respect to the Euclidean distance
on $\mathbb R^d$. That is, the $W_{ni}^\star(\mathbf{x^\star})$ are the
usual Euclidean NN weights \cite{gykokrha}, indexed by the random
set~$\mathcal L_n$.
\end{fact}

Recall that $|\mathcal R_n|$ represents the number of users who have
already provided information about the variable of interest (the movie
\textit{Titanic} in our example) at time~$n$. We are now in a position
to state the main result of this section.
\begin{theorem}
\label{convergence}
Suppose that $|{M}|\geq2$ and that assumption (\ref{eqF}) is
satisfied. Suppose that $k_n \to\infty$, $|\mathcal R_n|\to\infty$
a.s. and $\mathbb E[k_n/|\mathcal R_n|] \to0$ as $n\to\infty$. Then
\[
\mathbb E | \eta_n(\mathbf X^{\star})-\eta(\mathbf X^{\star}) |\to0
\qquad\mbox{as } n \to\infty.
\]
\end{theorem}

Thus, to achieve consistency, the number of nearest neighbors $k_n$,
over which one averages in order to estimate the regression function,
should on one hand, tend to infinity but should, on the other hand, be
small with respect to the cardinality of the subset of database users
who have already rated the item of interest. We illustrate this result
by working out two examples.
\begin{exmp}
\label{exemple1}
Consider, to start with, the somewhat ideal situation where all users
in the database have rated the item of interest. In this case,
$\mathcal R_n=\{1,\ldots,n\}$, and the asymptotic conditions on $k_n$
become $k_n\to\infty$ and $k_n/n\to0$ as $n\to\infty$. These are just
the well-known conditions ensuring consistency of the usual (i.e.,
Euclidean) NN regression estimate (\cite{gykokrha}, Chapter 6).
\end{exmp}
\begin{exmp}
\label{exemple2}
In this more sophisticated model, we recursively define the sequence
$(\mathcal{R}_n)_n$ as follows. Fix, for simplicity, $\mathcal{R}_1=\{
1\}$. At step $n\geq2$, we first decide (or not) to add one element to
$\mathcal{R}_{n-1}$ with probability $p\in(0,1)$, independently of
the data. If we decide to increase $\mathcal R_n$, then we do it by
picking a random variable $B_n$ uniformly over the set $\{1,\ldots,n\}
-\mathcal{R}_{n-1}$, and set $\mathcal R_n=\mathcal R_{n-1}\cup\{B_n\}
$; otherwise, $\mathcal R_n=\mathcal R_{n-1}$. Clearly, $|\mathcal
{R}_n|-1$ is a sum of $n-1$ independent Bernoulli random variables with
parameter $p$, and it has therefore a binomial distribution with
parameters $n-1$ and $p$. Consequently,
\[
\mathbb{E}\biggl[\frac{k_n}{|{\mathcal R_n}|} \biggr]=\frac{k_n
[1-(1-p)^n ]}{np}.
\]
In this setting, consistency holds provided $k_n \to\infty$ and
$k_n=\mbox{o}(n)$ as $n\to\infty$.
\end{exmp}

In the sequel, the letter $C$ will denote a positive constant, the
value of which may vary from line to line. Proof of Theorem \ref
{convergence} will strongly rely on Facts \ref{FAIT1}, \ref{FAIT2} and
the following proposition.
\begin{pro}
\label{pro:dec_eta}
Suppose that ${|{M}|}\geq2$ and that assumption (\ref{eqF}) is
satisfied. Let $\alpha_{ni}=\mathbb{P}(M^{n+1-i}\not\supset M | M)$. Then
\begin{eqnarray*}
&& \mathbb{E}|\eta_n(\mathbf{X}^\star)-\eta(\mathbf{X}^\star) |\\
&&\qquad \leq C \biggl\{\mathbb{E}\biggl[\frac{k_n}{|{\mathcal R_n}|} \biggr]
+\mathbb{E}\biggl[\frac{1}{|{\mathcal R_n}|}\sum_{i\in\mathcal
R_n}\mathbb{E}
\alpha_{ni} \biggr] +\mathbb{E}\biggl[\prod_{i\in\mathcal R_n}\alpha
_{ni} \biggr]\\
&&\hspace*{68.9pt}\qquad\quad{} + \mathbb{E}\biggl|\sum_{i\in\mathcal L_n}W_{ni}^\star(\mathbf{X}^\star
)\frac
{Y_i}{\|\mathbf{X}_i^\star\|}-\varphi(\mathbf{X}^\star) \biggr| \biggr\},
\end{eqnarray*}
where $\mathcal R_n$ stands for the nonempty subset of users who have
already provided information about the variable of interest at time
$n$, and $\mathcal L_n$ is defined in (\ref{eq:defLn}).
\end{pro}
\begin{pf}
Since $\|\mathbf{X}^\star\|\leq s\sqrt{d}$, it will be enough to upper
bound the quantity
\[
\mathbb{E}\biggl|\sum_{i\in\mathcal R_n}W_{ni}(\mathbf{X}^\star)\frac
{Y_i}{\|
\mathbf{X}_i^{(n)}\|}-\varphi(\mathbf{X}^\star)\biggr|.
\]
To this aim, we write
\begin{eqnarray*}
&& \mathbb{E}\biggl|\sum_{i\in\mathcal R_n}W_{ni}(\mathbf{X}^\star)\frac
{Y_i}{\|
\mathbf{X}_i^{(n)}\|}-\varphi(\mathbf{X}^\star) \biggr|\\
&&\qquad \leq\mathbb{E}\biggl[ \sum_{i\in\mathcal L_n^c}W_{ni}(\mathbf{X}^\star
)\frac
{Y_i}{\|\mathbf{X}_i^{(n)}\|} \biggr]+\mathbb{E}\biggl|\sum_{i\in\mathcal L
_n}W_{ni}(\mathbf{X}^\star)\frac{Y_i}{\|\mathbf{X}_i^{(n)}\|}-\varphi(\mathbf{X}^\star)\biggr|,
\end{eqnarray*}
where the symbol $A^c$ denotes the complement of the set $A$.
Let the event
\[
\mathcal A_n= \bigl[\exists i\in\mathcal L_n^c\dvtx\mathbf{X}_i^{(n)}\mbox
{ is among
the $k_n$-MS of }
\mathbf{X}^\star\mbox{ in }\bigl\{\mathbf{X}_i^{(n)}, i\in\mathcal R_n\bigr\} \bigr].
\]
Since $\sum_{i\in\mathcal L_n^c}W_{ni}(\mathbf{X}^\star)\leq1$, we have
\[
\mathbb{E}\biggl[\sum_{i\in\mathcal L_n^c}W_{ni}(\mathbf{X}^\star)\frac
{Y_i}{\|\mathbf{X}_i^{(n)}\|} \biggr] =\mathbb{E}\biggl[\sum_{i\in\mathcal L_n^c}W_{ni}(\mathbf{X}^\star
)\frac{Y_i}{\|\mathbf{X}_i^{(n)}\|}\mathbf{1}_{\mathcal A_n} \biggr]\leq
s\mathbb{P}
(\mathcal A_n).
\]
Observing that, for $i\in\mathcal L_n$, $\mathbf{X}_i^{(n)}=\mathbf{X}_i^\star$ and
$W_{ni}(\mathbf{X}^\star)\mathbf{1}_{\mathcal A_n^c}=W_{ni}^\star(\mathbf{X}^\star
)\mathbf{1}_{\mathcal A_n^c}$ (Fact \ref{FAIT2}), we obtain
\begin{eqnarray*}
&& \mathbb{E}\biggl|\sum_{i\in\mathcal L_n}W_{ni}(\mathbf{X}^\star)\frac
{Y_i}{\|\mathbf{X}_i^{(n)}\|}-\varphi(\mathbf{X}^\star) \biggr| \\
&&\qquad = \mathbb{E}\biggl|\sum_{i\in\mathcal L_n}W_{ni}(\mathbf{X}^\star)\frac
{Y_i}{\|\mathbf{X}_i^\star\|}-\varphi(\mathbf{X}^\star) \biggr| \\
&&\qquad = \mathbb{E}\biggl|\sum_{i\in\mathcal L_n}W_{ni}(\mathbf{X}^\star)\frac
{Y_i}{\|\mathbf{X}_i^\star\|}-\varphi(\mathbf{X}^\star) \biggr|\mathbf{1}_{\mathcal A_n} \\
&&\qquad\quad{} +\mathbb{E}\biggl|\sum_{i\in\mathcal L_n}W_{ni}^{\star}(\mathbf{X}^\star
)\frac
{Y_i}{\|\mathbf{X}_i^\star\|}-\varphi(\mathbf{X}^\star) \biggr|\mathbf
{1}_{\mathcal
A_n^c} \\
&&\qquad \leq s\mathbb{P}(\mathcal A_n)+\mathbb{E}\biggl|\sum_{i\in\mathcal
L_n}W_{ni}^{\star
}(\mathbf{X}^\star)\frac{Y_i}{\|\mathbf{X}_i^\star\|}-\varphi(\mathbf{X}^\star
) \biggr|.
\end{eqnarray*}
Applying finally Lemma \ref{lem:probAn} completes the proof of the proposition.
\end{pf}

We are now in a position to prove Theorem \ref{convergence}.
\begin{pf*}{Proof of Theorem \protect\ref{convergence}}
According to Proposition \ref{pro:dec_eta}, Lemma \ref{lem:Lnvide} and
Lemma \ref{lem:LnsurRn}, the result will be proven if we show that
\[
\mathbb{E}\biggl|\sum_{i\in\mathcal L_n}W_{ni}^\star(\mathbf{X}^\star)\frac
{Y_i}{\|\mathbf{X}_i^\star\|}-\varphi(\mathbf{X}^\star) \biggr|\to0 \qquad\mbox{as }
n\to\infty.
\]
For $L_n\in\mathcal P(\{1,\ldots,n\})$, set
\begin{eqnarray*}
Z_{L_n}^n &=& \frac{1}{k_n}\sum_{i\in L_n}\mathbf{1}_{ \{{\mathbf{X}_i^\star}/{\|\mathbf{X}_i^\star\|}\ \mathrm{is}
\ \mathrm{among}\ \mathrm{the}\ k_n\mbox{-}\mathrm{NN}\ \mathrm{of}\
{\mathbf{X}^\star}/{\|\mathbf{X}^\star\|}\ \mathrm{in}\ \{
{\mathbf{X}_i^\star}/{\|\mathbf{X}_i^\star\|}, i\in
L_n \} \}}\frac{Y_i}{\|\mathbf{X}_i^\star\|}\\
&&{} -\varphi(\mathbf{X}^\star).
\end{eqnarray*}
Conditionally on the event $[M=m]$, the random variables $\mathbf{X}^\star
$ and $\{\mathbf{X}_i^\star, i\in L_n\}$ are independent and identically
distributed. Thus, applying Theorem 6.1 in \cite{gykokrha}, we obtain
\[
\forall\varepsilon>0\qquad \exists A_m\geq1\dvtx k_n\geq A_m\quad
\mbox{and}\quad
\frac
{|{L_n}|}{k_n}\geq A_m\quad\Longrightarrow\quad\mathbb{E}_m|Z_{L_n}^n|\leq
\varepsilon,
\]
where we use the notation $\mathbb{E}_m[\cdot]=\mathbb{E}[\cdot|M=m]$. Let
$\mathbb{P}
_m(\cdot)=\mathbb{P}
(\cdot|M=m)$. By independence,
\[
\mathbb{E}_m|Z_{\mathcal L_n}^n| = \sum_{L_n\in\mathcal P(\{1,\ldots
,n\}
)}\mathbb{E}
_m|Z_{L_n}^n| \mathbb{P}_m(\mathcal L_n=L_n).
\]
Consequently, letting $A=\max A_m$, where the maximum is taken over all
possible choices of $m \in\mathcal P^{\star}(\{1,\ldots, d\})$, we
get, for all $n$ such that $k_n\geq A$,
\begin{eqnarray*}
\mathbb{E}_m|Z_{\mathcal L_n}^n| & = & \mathop{\sum_{L_n\in\mathcal
P(\{
1,\ldots
,n\})}}_{|{L_n}|\geq Ak_n} \mathbb{E}_m|Z_{L_n}^n| \mathbb
{P}_m(\mathcal L_n=L_n)\\
&&{} + \mathop{\sum_{L_n\in\mathcal P(\{1,\ldots,n\})}}_{|{L_n}|<
Ak_n} \mathbb{E}_m|Z_{L_n}^n| \mathbb{P}_m(\mathcal L_n=L_n) \\
&\leq& \varepsilon+s\mathbb{P}_m(|{\mathcal L_n}|<Ak_n).
\end{eqnarray*}
Therefore,
\[
\mathbb{E}|Z_{\mathcal L_n}^n|=\mathbb{E}[\mathbb{E}[ |Z_{\mathcal
L_n}^n| | M
] ] \leq\varepsilon+ s \mathbb{P}(|{\mathcal L_n}|<Ak_n ).
\]
Moreover, by Lemma \ref{lem:LnsurRn},
\[
\frac{|{\mathcal L_n}|}{k_n}=\frac{|{\mathcal R_n}|}{k_n} \biggl(1- \frac
{|{\mathcal L_n^c}|}
{|{{\mathcal R}_n}| } \biggr)\to\infty\qquad\mbox{in probability as }
n\to\infty.
\]
Thus for all $\varepsilon>0$, $\limsup_{n\to\infty}\mathbb
{E}|Z_{\mathcal L
_n}^n|\leq
\varepsilon$, whence $\mathbb{E}|Z_{\mathcal L_n}^n|\to0$ as $n\to
\infty$. This
shows the desired result.
\end{pf*}

%s5 ###
\section{Rates of convergence}\label{sec5}

In this section, we bound the rate of convergence of $\mathbb E
|\eta_n(\mathbf X^{\star})-\eta(\mathbf X^{\star}) |$ for the
cosine-type $k_n$-NN regression estimate. To reach this objective, we
will require that the function
\[
\varphi(\mathbf{x}^\star)= \mathbb E \biggl[ \frac{Y}{\|\mathbf{X}^\star\|} \bigg|
\frac{\mathbf{X}^\star}{\|\mathbf{X}^\star\|}=\frac{\mathbf{x}^\star}{\|
\mathbf{x}^\star\|} \biggr],
\]
satisfies a Lipschitz-type property with respect to the similarity
$\bar S$. More precisely, we say that $\varphi$ is Lipschitz with
respect to $\bar S$ if there exists a constant $C>0$ such that, for all
$\mathbf{x}$ and $\mathbf{x}'$ in $\mathbb R^d$,
\[
|\varphi(\mathbf{x})-\varphi(\mathbf{x}')|\leq C \sqrt{1-{\bar
S}(\mathbf{x},\mathbf{x}')}.
\]
In particular, for $\mathbf{x}$ and $\mathbf{x}' \in\mathbb R^d-\mathbf{0}$ with the
same null components, this property can be rewritten as
\[
|\varphi(\mathbf{x})-\varphi(\mathbf{x}')|\leq\frac{C}{\sqrt2}
\dd \biggl(
\frac{\mathbf{x}}{\|\mathbf{x}\|},\frac{\mathbf{x}'}{\|\mathbf
{x}'\|} \biggr),
\]
where we recall that $\dd$ denotes Euclidean distance.
\begin{theorem}
\label{vitesse}
Suppose that assumption (\ref{eqF}) is satisfied and that $\varphi$
is Lipschitz with respect to ${\bar S}$. Let $\alpha_{ni}=\mathbb{P}
(M^{n+1-i}\not\supset M | M)$, and assume that $|{M}|\geq4$. Then
there exists $C>0$ such that, for all $n \geq1$,
\begin{eqnarray*}
&& \mathbb E |\eta_n(\mathbf X^{\star})-\eta(\mathbf X^{\star})
|\\
&&\qquad \leq C \biggl\{\mathbb E \biggl[ \frac{k_n}{|\mathcal R_n|}\sum
_{i \in\mathcal R_n} \mathbb\mathbb{E}\alpha_{ni} \biggr] +\mathbb
{E}\biggl[\prod
_{i\in\mathcal R_n}\alpha_{ni} \biggr]+\mathbb E \biggl[ \biggl(\frac
{k_n}{|\mathcal R_n|} \biggr)^{P_n} \biggr]+\frac{1}{\sqrt{k_n}}
\biggr\},
\end{eqnarray*}
where $P_n=1/(|{M}|-1)$ if $k_n\leq|{\mathcal R_n}|$, and
$P_n=1$ otherwise.
\end{theorem}

To get an intuition on the meaning of Theorem \ref{vitesse}, it helps
to note that the terms depending on $\alpha_{ni}$ do measure the
influence of the unrated items on the performance of the estimate.
Clearly, this performance improves as the $\alpha_{ni}$ decrease, that
is, as the proportion of rated items growths. On the other hand, the
term $\mathbb E [ ({k_n}/{|\mathcal R_n|} )^{P_n}]$ can be interpreted
as a bias term in dimension $|M|-1$, whereas ${1}/{\sqrt{k_n}}$
represents a variance term. As usual in nonparametric estimation, the
rate of convergence of the estimate is dramatically deteriorated as
$|M|$ becomes large. However, in practice, this drawback may be
circumvented by using preliminary dimension reduction steps, such as
factorial methods (PCA, etc.) or inverse regression methods (SIR, etc.).
\begin{exmp}[(Example \ref{exemple1}, continued)]
Recall that we assume, in this ideal model, that $\mathcal R_n=\{
1,\ldots, n\}$. Suppose in addition that $M=\{1,\ldots,d\}$, that is,
any new user in the database rates all products the first time he
enters the database. Then the upper bound of Theorem \ref{vitesse} becomes
\[
\mathbb E |\eta_n(\mathbf X^{\star})-\eta(\mathbf X^{\star})
|= \mathrm{O} \biggl( \biggl(\frac{k_n}{n} \biggr)^{1/(d-1)}+\frac
{1}{\sqrt{k_n}} \biggr).
\]
Since neither $\mathcal R_n$ nor $M$ are random in this model, we see
that there is no influence of the dynamical rating process. Besides, we
recognize the usual rate of convergence of the Euclidean NN regression
estimate (\cite{gykokrha}, Chapter 6) in dimension $d-1$. In
particular, the choice $k_n\sim n^{2/(d+1)}$ leads to
\[
\mathbb E |\eta_n(\mathbf X^{\star})-\eta(\mathbf X^{\star})
|=\mathrm{O} \bigl(n^{-1/(d+1)} \bigr).
\]
Note that we are led to a ($d-1$)-dimensional rate of convergence
(instead of the usual $d$) just because everything happens as if the
data is projected on the unit sphere of $\mathbb{R}^d$.
\end{exmp}
\begin{exmp}[(Example \ref{exemple2}, continued)]
In addition to model \ref{exemple2}, we suppose that at each time, a
user entering the game reveals his preferences according to the
following sequential procedure. At time 1, the user rates exactly 4
items by randomly guessing in $\{1,\ldots,d\}$. At time $2$, he updates
his preferences by adding exactly one rating among his unrated items,
randomly chosen in $\{1,\ldots,d\}-M_1^1$. Similarly, at time $3$, the
user revises his preferences according to a new item uniformly selected
in $\{1,\ldots,d\}-M_1^2$, and so on. In such a scenario, $|M^j|=\min
(d,j+3)$ and thus, $M^j=\{1,\ldots,d\}$ for $j\geq d-3$. Moreover,
since $|M|=4$, a moment's thought shows that
\[
\alpha_{ni}=
\cases{
0, &\quad if $i\leq n-d+4$, \cr
1-\frac{ {d-4\choose n-i}}{ {d\choose n+4-i}},
&\quad if $n-d+5\leq i\leq n$.}
\]
Assuming $n\geq d-5$, we obtain
\begin{eqnarray*}
\sum_{i\in\mathcal R_n}\alpha_{ni} & \leq &\sum_{i=n-d+5}^n\alpha
_{ni}\\
&\leq&\sum_{i=n-d+5}^n \biggl(1-\frac
{(n+4-i)(n+3-i)(n+2-i)(n+1-i)}{d(d-1)(d-2)(d-3)} \biggr) \\
& \leq&(d-4) \biggl(1-\frac{24}{d(d-1)(d-2)(d-3)} \biggr).
\end{eqnarray*}
Similarly, letting $\mathcal{R}_{n0}=\mathcal R_n\cap\{n-d+5,\ldots
,n\}$, we have
\begin{eqnarray*}
\prod_{i\in\mathcal R_n}\alpha_{ni}&=& \prod_{i\in\mathcal
{R}_{n0}}\alpha_{ni} \mathbf{1}_{\{\min(\mathcal R_n)\geq n-d+5\}}\\
&\leq&\biggl(1-\frac{24}{d(d-1)(d-2)(d-3)} \biggr)^{|{\mathcal
{R}_{n0}}|}\mathbf{1}_{\{\min(\mathcal R_n)\geq n-d+5\}}.
\end{eqnarray*}
Since $|{\mathcal R}_n|-1$ has binomial distribution with parameters
$n-1$ and $p$, we obtain
\begin{eqnarray*}
\mathbb{E}\biggl[\prod_{i\in\mathcal R_n}\alpha_{ni} \biggr] &\leq& \mathbb
{P} \bigl(\min(\mathcal R_n)\geq n-d+5 \bigr)\\
& \leq& \mathbb{P} (|\mathcal{R}_n|\leq d-5 ) \leq\frac{C}{n}.
\end{eqnarray*}
Finally, applying Jensen's inequality,
\begin{eqnarray*}
\mathbb E \biggl[ \biggl(\frac{k_n}{|\mathcal R_n|} \biggr)^{P_n} \biggr]
& = & \mathbb E \biggl[ \biggl(\frac{k_n}{|\mathcal R_n|} \biggr)^{1/3}
\mathbf{1}_{\{k_n\leq|\mathcal{R}_n|\}} \biggr]+\mathbb E\biggl[ \frac
{k_n}{|\mathcal R_n|}\mathbf{1}_{\{k_n > |\mathcal{R}_n|\}} \biggr]\\
& \leq& C \biggl(\mathbb E \biggl[ \frac{k_n}{|\mathcal R_n|} \biggr]
\biggr)^{1/3}\leq C \biggl(\frac{k_n}{n} \biggr)^{1/3}.
\end{eqnarray*}
Putting all the pieces together, we get with Theorem \ref{vitesse}
\[
\mathbb E |\eta_n(\mathbf X^{\star})-\eta(\mathbf X^{\star})
| = \mathrm{O} \biggl( \biggl(\frac{k_n}{n} \biggr)^{1/3}+\frac
{1}{\sqrt{k_n}} \biggr).
\]
In particular, the choice $k_n \sim n^{2/5}$ leads to
\[
\mathbb E |\eta_n(\mathbf X^{\star})-\eta(\mathbf X^{\star})
| =\mathrm{O}(n ^{-1/5}),
\]
which is the usual NN regression estimate rate of convergence when the
data is projected on the unit sphere of $\mathbb{R}^4$.
\end{exmp}
\begin{pf*}{Proof of Theorem \protect\ref{vitesse}}
Starting from Proposition \ref{pro:dec_eta}, we just need to upper
bound the quantity
\[
\mathbb{E}\biggl|\sum_{i\in\mathcal L_n}W_{ni}^\star(\mathbf{X}^\star)\frac
{Y_i}{\|\mathbf{X}_i^\star\|}-\varphi(\mathbf{X}^\star) \biggr|.
\]
A combination of Lemma \ref{lem:rate1ppv} and the proof of Theorem 6.2
in \cite{gykokrha} shows that
%
%e5.1 ###
%
\begin{eqnarray}\label{ineq:M}
&& \mathbb{E}\biggl|\sum_{i\in\mathcal L_n}W_{ni}^\star(\mathbf{X}^\star
)\frac{Y_i}{\|
\mathbf{X}_i^\star\|}-\varphi(\mathbf{X}^\star)
\biggr|\nonumber\\[-8pt]\\[-8pt]
&&\qquad \leq C \biggl\{\frac{1}{\sqrt{k_n}}+\mathbb{E}\biggl[ \biggl(\frac
{k_n}{|{\mathcal L_n}|} \biggr)^{1/(|{M}|-1)}\mathbf{1}_{\{\mathcal L_n\neq
\varnothing
\}} \biggr]+\mathbb{P}(\mathcal L_n=\varnothing) \biggr\}.\nonumber
\end{eqnarray}
We obtain
\begin{eqnarray*}
&& \mathbb{E}\biggl[ \biggl(\frac{k_n}{|{\mathcal L_n}|} \biggr)^{1/(|{M}|-1)}
\mathbf{1}_{\{\mathcal L_n\neq\varnothing\}} \biggr]\\
&&\qquad =\mathbb{E}\biggl[ \biggl(\frac{k_n}{{|{\mathcal{R}_n}|} (1-|{\mathcal L
_n^c}|/|{\mathcal{R}_n}| )} \biggr)^{1/({|{M}|}-1)}\mathbf{1}_{\{|{\mathcal
L_n^c}|\leq|{\mathcal{R}_n}|/2\}} \biggr]\\
&&\qquad\quad{} + \mathbb{E}\biggl[ \biggl(\frac{k_n}{{|{\mathcal L_n}|}}
\biggr)^{1/({|{M}|}-1)}\mathbf{1}_{\{|{\mathcal L_n^c}|> |{\mathcal {R}_n}|/2\}}
\mathbf{1}_{\{\mathcal L_n\neq\varnothing\}} \biggr]\\
&&\qquad \leq\mathbb{E}\biggl[ \biggl( \frac{2k_n}{{|{\mathcal R_n}|}} \biggr)^{1/(|{M}|-1)}
\biggr]+\mathbb{E}\bigl[ k_n^{1/(|{M}|-1)}
\mathbf{1}_{\{|{\mathcal L_n^c}|>|{\mathcal R_n}|/2\}} \bigr].
\end{eqnarray*}
Since $|{M}|\geq4$, one has $2^{1/({|{M}|}-1)}\leq2$ and
$k_n^{1/(|{M}|-1)} \leq k_n$ in the rightmost term, so that, thanks
to Lemma
\ref{lem:LnsurRn},
\begin{eqnarray*}
&& \mathbb{E}\biggl[ \biggl(\frac{k_n}{|{\mathcal L_n}|} \biggr)^{1/(|{M}|-1)}
\mathbf{1}_{\{\mathcal L_n\neq\varnothing\}} \biggr]\\
&&\qquad \leq C
\biggl\{ \mathbb{E}\biggl[ \biggl( \frac{k_n}{{|{\mathcal R_n}|}}
\biggr)^{1/(|{M}|-1)} \biggr]+\mathbb{E}\biggl[\frac{k_n}{|{\mathcal R_n}|}\sum
_{i\in\mathcal R_n}\mathbb{E}\alpha_{ni} \biggr] \biggr\}.
\end{eqnarray*}
The theorem is a straightforward combination of Proposition \ref
{pro:dec_eta}, inequality (\ref{ineq:M}) and Lemma \ref{lem:Lnvide}.
\end{pf*}

%s6 ###
\section{Technical lemmas}\label{sec6}

Before stating some technical lemmas, we remind the reader that
$\mathcal R_n$ stands for the nonempty subset of $\{1,\ldots,n\}$ of
users who have already rated the variable of interest at time $n$.
Recall also that, for all $i\geq1$,
\[
T_i=\min(k\geq i\dvtx M_i^{k+1-i}\supset M)
\]
and
\[
\mathcal L_n=\{i\in\mathcal R_n\dvtx T_i\leq n\}.
\]
\begin{lem}
\label{lem:Lnvide}
We have
\[
\mathbb{P}(\mathcal L_n=\varnothing)=\mathbb{E}\biggl[\prod_{i\in
\mathcal R_n}\alpha
_{ni} \biggr]\to0 \qquad\mbox{as } n\to\infty.
\]
\end{lem}
\begin{pf}
Conditionally on $M$ and $\mathcal{R}_n$, the random variables $\{
T_i,i\in
\mathcal{R}_n\}$ are independent. Moreover, the sequence $(M^n)_{n\geq1}$
is nondecreasing. Thus, the identity $[T_i>n]=[M_i^{n+1-i}\not\supset
M]$ holds for all $i\in\mathcal R_n$. Hence,
\begin{eqnarray}
\mathbb{P}(\mathcal L_n=\varnothing) & = &\mathbb{P}(\forall i\in
\mathcal
R_n\dvtx T_i>n )\nonumber\\
& = &\mathbb{E}[\mathbb{P}(\forall i\in\mathcal R_n\dvtx T_i>n |
\mathcal R_n,M ) ] \nonumber\\
& = &\mathbb{E}\biggl[\prod_{i\in\mathcal R_n}\mathbb{P}(T_i>n |
\mathcal R_n,M ) \biggr]\nonumber\\
& = &\mathbb{E}\biggl[\prod_{i\in\mathcal R_n}\mathbb{P}(M_i^{n+1-i}\not
\supset
M | M ) \biggr]\nonumber \\
\eqntext{\mbox{[by independence of $(M_i^{n+1-i},M)$ and $\mathcal
R_n$]}}\\
& = & \mathbb{E}\biggl[\prod_{i\in\mathcal R_n}\alpha_{ni} \biggr].\nonumber
\end{eqnarray}
The last statement of the lemma is clear since, for all $i$, $\alpha
_{ni}\to0$ a.s. as $n\to\infty$.
%\rightqed
\end{pf}
\begin{lem}
\label{lem:LnsurRn} We have
\[
\mathbb{E}\biggl[\frac{|{\mathcal L_n^c}|}{|{\mathcal R_n}|} \biggr]=\mathbb{E}
\biggl[\frac{1}{|{\mathcal R_n}|}\sum_{i\in\mathcal R_n}\mathbb{E}\alpha
_{ni} \biggr]
\]
and
\[
\mathbb{E}\biggl[ \frac{1}{|{\mathcal L_n}|}\mathbf{1}_{\{\mathcal L_n\neq
\varnothing\}}
\biggr]\leq2\mathbb{E}\biggl[\frac{1}{|{\mathcal{R}_n}|} \biggr]+2\mathbb{E}\biggl[\frac
{1}{|{\mathcal R_n}|}\sum_{i\in\mathcal R_n}\mathbb{E}\alpha_{ni} \biggr].
\]

Moreover, if $\lim_{n\to\infty}|{\mathcal R_n}|=\infty$ a.s., then
\[
\lim_{n\to\infty}\mathbb{E}\biggl[\frac{|{\mathcal L_n^c}|}{|{\mathcal
R_n}|} \biggr]= 0.
\]
\end{lem}
\begin{pf}
First, using the fact that the sequence $(M^n)_{n\geq1}$ is
nondecreasing, we see that for all $i\in\mathcal R_n$,
$[T_i>n]=[M_i^{n+1-i}\not\supset M]$. Next, recalling that $\mathcal
R_n$ is independent of $T_i$ for fixed $i$, we obtain
\[
\mathbb{E}\biggl[\frac{|{\mathcal L_n^c}|}{|{\mathcal R_n}|} \bigg| \mathcal
R_n \biggr]=\frac{1}{|{\mathcal R_n}|}\mathbb{E}\biggl[\sum_{i\in\mathcal
R_n}\mathbf{1}_{\{T_i>n\}} \Big| \mathcal R_n \biggr]=\frac{1}{|{\mathcal
R_n}|}\sum_{i\in\mathcal R_n}\mathbb{P}(M_i^{n+1-i}\not\supset M)
\]
and this proves the first statement of the lemma. Now define $\mathcal
J_n=\{n+1-i,i\in\mathcal R_n\}$ and observe that
\[
\mathbb{E}\biggl[ \frac{|{\mathcal L_n^c}|}{|{\mathcal R_n}|} \biggr]=\mathbb{E}
\biggl[\frac{1}{|{\mathcal J_n}|}\sum_{j\in\mathcal J_n}\mathbb
{P}(M^j\not
\supset M) \biggr],
\]
where we used $|\mathcal J_n|=|\mathcal R_n|$. Since, by assumption,
$|{\mathcal J_n}|=|{\mathcal R_n}|\to\infty$ a.s. as $n \to
\infty
$ and $\mathbb{P}(M^j\not\supset M)\to0$ as $j\to\infty$, we obtain
\[
\lim_{n\to\infty}\frac{1}{|{\mathcal J_n}|}\sum_{j\in\mathcal
J_n}\mathbb{P}(M^j\not\supset M)=0 \qquad\mbox{a.s.}
\]
The conclusion follows by applying Lebesgue's dominated convergence
theorem. The second statement of the lemma is obtained from the
following chain of inequalities:
\begin{eqnarray*}%\label{eq:Ln}
\mathbb{E}\biggl[ \frac{1}{|{\mathcal L_n}|}\mathbf{1}_{\{\mathcal L_n\neq
\varnothing\}} \biggr]
& = & \mathbb{E}\biggl[ \frac{1}{{|{\mathcal{R}_n}|} (1-|{\mathcal L
_n^c}|/|{\mathcal{R}_n}| )}\mathbf{1}_{\{\mathcal L_n\neq\varnothing\}}
\biggr]\\
& = & \mathbb{E}\biggl[ \frac{1}{{|{\mathcal{R}_n}|} (1-|{\mathcal L
_n^c}|/|{\mathcal{R}_n}| )}\mathbf{1}_{\{|{\mathcal L_n^c}|\leq
|{\mathcal {R}_n}|/2\}
} \biggr]\\
& &{} + \mathbb{E}\biggl[ \frac{1}{{|{\mathcal L_n}|}}\mathbf{1}_{\{|{\mathcal L_n^c}|>
|{\mathcal{R}_n}|/2\}} \mathbf{1}_{\{\mathcal L_n\neq\varnothing\}} \biggr]\\
& \leq& 2\mathbb{E}\biggl[ \frac{1}{|{\mathcal{R}_n}|} \biggr]+\mathbb{P}
\biggl(|{\mathcal L_n^c}|>\frac{|{\mathcal{R}_n}|}{2} \biggr)\\
& \leq& 2\mathbb{E}\biggl[ \frac{1}{|{\mathcal{R}_n}|} \biggr]+2\mathbb{E}
\biggl[\frac{|{\mathcal L_n^c}|}{|{\mathcal R_n}|} \biggr].
\end{eqnarray*}
Applying the first part of the lemma completes the proof.
\end{pf}
\begin{lem}
\label{lem:beurk} Denote by $\mathbf Z^\star$ and $\mathbf Z_1^\star$
the random variables $\mathbf Z^\star=\mathbf{X}^\star/\|\mathbf
{X}^\star\|$,
$\mathbf Z_1^\star=\mathbf{X}_1^\star/\|\mathbf{X}_1^\star\|$, and
let $\xi
(\mathbf
Z^{\star})=\mathbb{P}(S(\mathbf Z^{\star},\mathbf Z_1^\star)>1/2 |
\mathbf
Z^{\star})$. Then
\begin{eqnarray*}
\mathbb{P}\bigl(2k_n>|{\mathcal L_n}|\xi(\mathbf Z^\star) | \mathcal
L_n,M \bigr)
&\leq&
2\mathbb{E}\biggl[\frac{k_n}{|{\mathcal{R}_n}|} \bigg| \mathcal L_n \biggr] \mathbb{E}
\biggl[\frac{1}{\xi(\mathbf Z^\star)} \bigg| M \biggr] \\
&&{} +\mathbb{E}\biggl[\frac{|{\mathcal L_n^c}|}{|{\mathcal{R}_n}|} \bigg|
\mathcal L
_n,M \biggr].
\end{eqnarray*}
\end{lem}
\begin{pf}
If $M$ is fixed, $\mathbf Z^\star$ is independent of $\mathcal L_n$ and
$\mathcal R_n$. Thus by Markov's inequality,
\begin{eqnarray*}
&& \mathbb{P}\bigl(2k_n>|{\mathcal L_n}|\xi(\mathbf Z^\star) | \mathcal
L_n,M,\mathcal
{R}_n \bigr) \\
&&\qquad = \mathbb{P}\bigl(2k_n> |{\mathcal{R}_n}|\xi(\mathbf Z^\star
)-|{\mathcal L_n^c}|\xi(\mathbf Z^\star) | \mathcal L_n,M,\mathcal
{R}_n \bigr)\\
&&\qquad = \mathbb{P}\bigl(2k_n+ |{\mathcal L_n^c}|\xi(\mathbf Z^\star)\geq
|{\mathcal{R}_n}| \xi(\mathbf Z^\star) | \mathcal L_n,M,\mathcal
{R}_n \bigr)\\
&&\qquad \leq\frac{2k_n}{|{\mathcal{R}_n}|} \mathbb{E}\biggl[ \frac{1}{\xi
(\mathbf Z^\star)} \bigg| M \biggr]+\frac{|{\mathcal L_n^c}|}{|{\mathcal{R}_n}|}.
\end{eqnarray*}
The proof is completed by observing that $\mathcal{R}_n$ and $M$ are
independent random variables.
\end{pf}

Let $\mathcal B(\mathbf{x},\varepsilon)$ be the closed Euclidean ball in
$\mathbb R^d$ centered at $\mathbf{x}$ of radius $\varepsilon$. Recall
that the
support of a probability measure $\mu$ is defined as the closure of the
collection of all $\mathbf{x}$ with $\mu(\mathcal B(\mathbf{x},\varepsilon))>0$
for all $\varepsilon>0$. The next lemma can be proved with a slight
modification of the proof of Lemma 10.2 in \cite{DGL}.
\begin{lem}
\label{lem:intunsurmu}
Let $\mu$ be a probability measure on $\mathbb R^d$ with a compact
support. Then
\[
\int\frac{1}{\mu(\mathcal B(\mathbf{x},r))}\mu(\dd\mathbf{x})\leq C
\]
with $C>0$ a constant depending upon $d$ and $r$ only.
\end{lem}
\begin{lem}
\label{lem:probAn}
Suppose that $|{M}|\geq2$, and let the event
\[
\mathcal A_n= \bigl[\exists i\in\mathcal L_n^c\dvtx\mathbf{X}_i^{(n)}\mbox
{ is among
the $k_n$-MS of }\mathbf{X}^\star\mbox{ in }\bigl\{\mathbf{X}_i^{(n)},i\in
\mathcal R_n\bigr\} \bigr].
\]
Then
\[
\mathbb{P}(\mathcal A_n)\leq C \biggl\{\mathbb E \biggl[\frac{k_n}{|{\mathcal
R_n}|} \biggr]+\mathbb{E}\biggl[\frac{1}{|{\mathcal R_n}|}\sum
_{i\in\mathcal R_n}\mathbb{E}\alpha_{ni} \biggr]+\mathbb{E}\biggl[\prod_{i\in
\mathcal R_n} \alpha_{ni} \biggr] \biggr\}.
\]
\end{lem}
\begin{pf}
Recall that, for a fixed $i \in\mathcal R_n$, the random variable
$\mathbf{X}
_i^{\star}=(X_{i1}^{\star}$,\break$\ldots,X_{id}^{\star})$ is defined by
\[
X_{ij}^{\star}=
\cases{
X_{ij}, &\quad if $j \in M$,\cr
0, &\quad otherwise,}
\]
and $\mathbf{X}_i^{(n)}= \mathbf{X}_i^{\star}$ as soon as $M\subset
M_{i}^{n+1-i}$.

We first prove the inclusion
%
%e6.1 ###
%
\begin{equation}
\label{eq:part1_lem4}
\mathcal A_n\subset[|{\{j\in\mathcal L_n\dvtx S(\mathbf{X}^\star
,\mathbf{X} _j^\star)>1/2\} }|\leq k_n ].
\end{equation}
Take $i\in\mathcal L_n^c$ such that $\mathbf{X}_i^{(n)}$ is among the
$k_n$-MS of
$\mathbf{X}
^\star$ in $\{\mathbf{X}_i^{(n)},i\in\mathcal R_n\}$. Then, for all
$j\in
\mathcal L
_n$ such that $S(\mathbf{X}^\star,\mathbf{X}_j^\star)>1/2$, we have
\[
S(\mathbf{X}^\star,\mathbf{X}_j^\star)>\tfrac{1}{2}\geq
p_i^{(n)}\bar S\bigl(\mathbf{X}
^\star,\mathbf{X}
_i^{(n)}\bigr)=S\bigl(\mathbf{X}^\star,\mathbf{X}_i^{(n)}\bigr)
\]
since $p_i^{(n)}\leq1-1/|{M}|\leq1/2$ if $|{M}|\geq2$. If
\[
|{\{j\in\mathcal L_n\dvtx S(\mathbf{X}^\star,\mathbf{X}_j^\star
)>1/2\}}|> k_n,
\]
then $\mathbf{X}_i^{(n)}$ is not among the $k_n$-MS of $\mathbf
{X}^\star$ among the
$\{\mathbf{X}_i^{(n)},i\in\mathcal R_n\}$. This contradicts the
assumption on
$\mathbf{X}_i^{(n)}$ and proves inclusion (\ref{eq:part1_lem4}).

Next, define $\mathbf Z^\star=\mathbf{X}^\star/\|\mathbf{X}^\star\|
$, $\mathbf
Z_i^\star=\mathbf{X}_i^\star/\|\mathbf{X}_i^\star\|$, $i=1,\ldots
,n$, and let
$\xi
(\mathbf Z^{\star})=\mathbb{P}(S(\mathbf Z^{\star},\mathbf Z_1^\star)>1/2
| \mathbf Z^\star)$. If $k_n-|{\mathcal L_n}|\xi(\mathbf{Z}^\star
)\leq
-(1/2)|{\mathcal L_n}|\xi(\mathbf{Z}^\star)$ and $\mathcal L_n\neq
\varnothing$, we deduce from
(\ref{eq:part1_lem4}) that
\begin{eqnarray*}
&& \mathbb P (\mathcal A_n | \mathcal L_n,\mathbf{Z}^\star)\\
&&\qquad \leq\mathbb{P}\biggl(\sum_{j\in\mathcal L_n}\mathbf{1}_{\{S(\mathbf
{Z}^{\star},\mathbf{Z}
_j^{\star})>1/2\}}\leq k_n \Big| \mathcal L_n,\mathbf{Z}^\star\biggr) \\
&&\qquad = \mathbb{P}\biggl(\sum_{j\in\mathcal L_n} \bigl(\mathbf{1}_{\{S(\mathbf
{Z}^\star,\mathbf{Z}_j^\star
)>1/2\}}-\xi(\mathbf{Z}^\star) \bigr)\leq k_n-|{\mathcal L_n}|\xi
(\mathbf{Z}^\star)
\Big| \mathcal L_n,\mathbf{Z}^\star\biggr) \\
&&\qquad \leq\mathbb{P}\biggl(\sum_{j\in\mathcal L_n} \bigl(\mathbf{1}_{\{S(\mathbf
{Z}^\star,\mathbf{Z}
_j^\star)>1/2\}}-\xi(\mathbf{Z}^\star) \bigr)\leq-\frac
{1}{2}|{\mathcal L_n}|\xi
(\mathbf{Z}^\star) \Big| \mathcal L_n,\mathbf{Z}^\star\biggr) \\
&&\qquad \leq\frac{4|{\mathcal L_n}|\xi(\mathbf{Z}^\star)}{ (|{\mathcal
L_n}|\xi(\mathbf{Z}
^\star) )^2}=\frac{4}{|{\mathcal L_n}|\xi(\mathbf{Z}^\star)}\qquad\mbox{(by Chebyshev's inequality).}
\end{eqnarray*}
In the last inequality, we use the fact that, since $\sigma(M)\subset
\sigma(\mathbf Z^\star)$, the random variables $\{\mathbf Z^\star
_i,i\in\mathcal L_n\}$ are independent conditionally on $\mathbf
Z^\star$ and
$\mathcal L_n$. Using again the inclusion $\sigma(M)\subset\sigma
(\mathbf
Z^\star)$, we obtain, on the event $[\mathcal L_n \neq0]$,
\begin{eqnarray*}
&& \mathbb P (\mathcal A_n | \mathcal L_n,M )\\
&&\qquad = \mathbb{E}[ \mathbb{P}( \mathcal A_n | \mathcal L_n,\mathbf
Z^\star) |\mathcal L_n,M ]\\
&&\qquad \leq\frac{4}{|{\mathcal L_n}|}\mathbb{E}\biggl[\frac{1}{\xi(\mathbf
{Z}^\star)}
\bigg| \mathcal L_n,M \biggr]+\mathbb{P}\biggl(k_n-|{\mathcal L_n}|\xi(\mathbf
{Z}^\star)>-\frac
{1}{2}|{\mathcal L_n}|\xi(\mathbf{Z}^\star) \Big| \mathcal L_n,M \biggr) \\
&&\qquad = \frac{4}{|{\mathcal L_n}|}\mathbb{E}\biggl[\frac{1}{\xi(\mathbf
{Z}^\star)}
\bigg| M \biggr]+\mathbb{P}\bigl(|{\mathcal L_n}|\xi(\mathbf{Z}^\star)<2k_n |
\mathcal L
_n,M \bigr).
\end{eqnarray*}
Applying Lemma \ref{lem:beurk}, on the event $[\mathcal L_n\neq
\varnothing]$,
\begin{eqnarray*}
\mathbb P (\mathcal A_n | \mathcal L_n,M )
&\leq&\frac{4}{|{\mathcal L_n}|}\mathbb{E}\biggl[\frac{1}{\xi(\mathbf
{Z}^\star)} \bigg|
M \biggr]+
2\mathbb{E}\biggl[\frac{k_n}{|{\mathcal{R}_n}|} \bigg| \mathcal L_n \biggr] \mathbb{E}
\biggl[\frac{1}{\xi(\mathbf Z^\star)} \bigg| M \biggr]\\
&&{}+\mathbb{E}\biggl[\frac{|{\mathcal
L_n^c}|}{|{\mathcal{R}_n}|} \bigg| \mathcal L_n,M \biggr].
\end{eqnarray*}
Moreover, by Fact \ref{FAIT1},
\[
\xi(\mathbf{Z}^\star)=\mathbb{P}\bigl(S(\mathbf{Z}^\star,\mathbf
{Z}_1^\star)>\tfrac{1}{2} \big|
\mathbf{Z}^\star\bigr)\geq\mathbb{P}\bigl(d^2(\mathbf{Z}^\star,\mathbf
{Z}_1^\star)\leq\tfrac
{1}{2} \big| \mathbf{Z}^\star\bigr).
\]
Thus, denoting by $\nu^M$ the distribution of $\mathbf{Z}^\star$
conditionally
to $M$, we deduce from Lemma \ref{lem:intunsurmu} that
\[
\mathbb{E}\biggl[\frac{1}{\xi(\mathbf{Z}^\star)} \bigg| M \biggr]\leq\int\frac
{1}{\nu^M(\mathcal B(\mathbf{z},1/\sqrt2))}\nu^M(\dd
\mathbf{z})\leq C,
\]
where the constant $C$ does not depend on $M$. Putting all the pieces
together, we obtain
\[
\mathbb{P}(\mathcal A_n)\leq C \biggl\{ \mathbb{E}\biggl[\frac{1}{|{\mathcal
L_n}|}\mathbf{1}_{\{
\mathcal L_n\neq\varnothing\}} \biggr]+\mathbb{E}\biggl[\frac{k_n}{|{\mathcal{R}_n}|}
\biggr]+\mathbb{E}\biggl[ \frac{|{\mathcal L_n^c}|}{|{\mathcal{R}_n}|}
\biggr] \biggr\}+\mathbb{P}(\mathcal L_n=\varnothing).
\]
We conclude the proof with Lemmas \ref{lem:Lnvide} and \ref{lem:LnsurRn}.
\end{pf}

In the sequel, we let $\mathbf{X}_{(1)}^\star,\ldots,
\mathbf{X}_{({|{\mathcal L _n}|})}^\star$ be the sequence $\{\mathbf
{X}_i^\star,i\in\mathcal L_n\}$ reordered according to decreasing
similarities $S(\mathbf{X}^{\star},\mathbf{X}_{i}^\star ),i\in\mathcal
L_n$, that is,
\[
S \bigl(\mathbf{X}^{\star},\mathbf{X}_{(1)}^\star\bigr)\geq\cdots\geq S
\bigl(\mathbf{X}
^{\star},\mathbf{X}_{(|{\mathcal L_n}|)}^\star\bigr).
\]
Lemma \ref{lem:rate1ppv} below states the rate of convergence to 1 of
$S(\mathbf{X}^{\star},\mathbf{X}_{(1)}^\star)$.
\begin{lem}
\label{lem:rate1ppv}
Suppose that $|{M}|\geq4$. Then there exists $C>0$ such that, on
the event $[\mathcal L_n\neq\varnothing]$,
\[
1-\mathbb{E}\bigl[S\bigl(\mathbf{X}^\star,\mathbf{X}_{(1)}^\star\bigr) | M,\mathcal L_n
\bigr]\leq
\frac{C}{|{\mathcal L_n}|^{2/(|{M}|-1)}}.
\]
\end{lem}
\begin{pf}
Observe that
\begin{eqnarray*}
&&\mathbb{E} \bigl[ 1-S\bigl(\mathbf{X}^\star,\mathbf{X}_{(1)}^\star\bigr) |
\mathbf{X}^\star,\mathcal L_n \bigr]
\\
&&\qquad = \int_0^1 \mathbb{P}\bigl(1-S\bigl(\mathbf{X}^\star,\mathbf
{X}_{(1)}^\star
\bigr)>\varepsilon| \mathbf{X}^\star,\mathcal L_n \bigr)\, \dd\varepsilon\\
&&\qquad =\int_0^1 \mathbb{P}\bigl( \forall i\in\mathcal L_n \dvtx1-S(\mathbf
{X}^\star
,\mathbf{X}_i^\star)>\varepsilon| \mathbf{X}^\star,\mathcal L_n
\bigr)\, \dd\varepsilon.
\end{eqnarray*}
Since $\sigma(M)\subset\sigma(\mathbf{X}^\star)$, given $\mathbf
{X}^\star$ and
$\mathcal L
_n$, the random variables $\{\mathbf{X}_i^\star,i\in\mathcal L_n\}$
are independent
and identically distributed. Hence,
\[
\mathbb{E}\bigl[ 1-S\bigl(\mathbf{X}^\star,\mathbf{X}_{(1)}^\star\bigr) | \mathbf
{X}^\star, \mathcal L_n \bigr]=
\int_0^1 \bigl[\mathbb{P}\bigl( 1-S(\mathbf{X}^\star,\mathbf{X}_1^\star
)>\varepsilon
| \mathbf{X}^\star\bigr) \bigr]^{|{\mathcal L_n}|}\,\dd\varepsilon.
\]
Denote by $\nu^M$ the conditional distribution of $\mathbf{X}^\star
/\|\mathbf{X}
^\star
\|$ given $M$. The support of $\nu^M$ is contained in both the unit
sphere of $\mathbb R^d$ and in a $|{M}|$-dimen\-sional vector space. Thus,
for simplicity, we shall consider that the support of $\nu^M$ is
contained in the unit sphere of $\mathbb R^{|{M}|}$. Let $\mathcal
{B}^{|{M}|}(\mathbf{x},r)$ be the closed Euclidean ball in $\mathbb R^{|{M}|}$
centered at
$\mathbf{x}$ of radius $r$. Since $\mathbf{X}^\star$ (resp.,
$\mathbf{X}^\star
_1$) only
depends on $M$ and $\mathbf{X}$ (resp., $\mathbf{X}_1$), then, given
$\mathbf{X}^\star$, the
random variable $\mathbf{X}^\star_1/\|\mathbf{X}_1^\star\|$ is
distributed according
to $\nu^M$. Thus, for any $\varepsilon>0$, we may write (Fact \ref{FAIT1})
\[
\mathbb{P}\bigl( 1-S(\mathbf{X}^\star,\mathbf{X}_1^\star)>\varepsilon|
\mathbf{X}^\star
\bigr)=1-\nu^M \biggl(\mathcal B^{|{M}|} \biggl(\frac{\mathbf{X}^{\star}}{\|\mathbf{X}
^{\star}\|},\sqrt{2\varepsilon} \biggr) \biggr)
\]
and, consequently,
\[
\mathbb{E}\bigl[ 1-S\bigl(\mathbf{X}^\star,\mathbf{X}_{(1)}^\star\bigr) | \mathbf
{X}^\star,\mathcal L_n \bigr]=
\int_0^1 \biggl[ 1-\nu^M \biggl(\mathcal B^{|{M}|} \biggl(\frac{\mathbf{X}^{\star
}}{\|\mathbf{X}^{\star}\|},\sqrt{2\varepsilon} \biggr) \biggr) \biggr]^{|{\mathcal L
_n}|}\,\dd\varepsilon.
\]
Using the inclusion $\sigma(M)\subset\sigma(\mathbf{X}^\star)$, we obtain
%
%e6.2 ###
%
\begin{eqnarray}
\label{eq:bla}
&& \mathbb{E}\bigl[ 1-S\bigl(\mathbf{X}^\star,\mathbf{X}_{(1)}^\star\bigr) |
M,\mathcal L_n \bigr]
\nonumber\\[-8pt]\\[-8pt]
&&\qquad = \int_0^1\mathbb{E}\biggl[ \biggl\{1-\nu^M \biggl(\mathcal B^{|{M}|} \biggl(\frac
{\mathbf{X}^{\star}}{\|\mathbf{X}^{\star}\|},\sqrt{2\varepsilon}
\biggr) \biggr) \biggr\}^{|{\mathcal L_n}|} \bigg| M,\mathcal L_n
\biggr]\,\dd\varepsilon.\nonumber
\end{eqnarray}

Fix $\varepsilon>0$, and denote by $\mathcal{S}(M)$ the support of
$\nu
^M$. There exists Euclidean balls $A_1,\ldots,A_{N(\varepsilon)}$ in
$\mathbb R
^{|{M}|}$ with radius $\sqrt{2\varepsilon}/2$ such that
\[
\mathcal{S}(M)\subset\bigcup_{j=1}^{N(\varepsilon)} A_j
\quad\mbox{and}\quad
N(\varepsilon)\leq\frac{C}{\varepsilon^{({|{M}|}-1)/2}}
\]
for some $C>0$ which may be chosen independently of $M$. Clearly, if
$\mathbf{x}\in A_j\cap\mathcal{S}(M)$, then $A_j\subset\mathcal
B^{|{M}|}(\mathbf{x},\sqrt{2\varepsilon})$. Thus
\begin{eqnarray*}
&& \mathbb{E}\biggl[ \biggl\{1-\nu^M \biggl(\mathcal B^{|{M}|} \biggl(\frac{\mathbf{X}
^{\star}}{\|\mathbf{X}^{\star}\|},\sqrt{2\varepsilon} \biggr) \biggr) \biggr\}
^{|{\mathcal L_n}|} \bigg|M,\mathcal L_n \biggr] \\
&&\qquad \leq\sum_{j=1}^{N(\varepsilon)}\int_{A_j}\mathbb{E}\biggl[ \biggl\{1-\nu
^M \biggl(\mathcal B^M \biggl(\frac{\mathbf{X}^{\star}}{\|\mathbf{X}^{\star}\|
},\sqrt
{2\varepsilon} \biggr) \biggr) \biggr\}^{|{\mathcal L_n}|} \bigg|M,\mathcal L_n \biggr]
\nu^M(\dd\mathbf{x}) \\
&&\qquad \leq\sum_{j=1}^{N(\varepsilon)}\int_{A_j} \bigl(1-\nu
^M(A_j) \bigr)^{|{\mathcal L_n}|}\nu^M(\dd\mathbf{x}) \\
&&\qquad \leq\sum_{j=1}^{N(\varepsilon)}\nu^M(A_j) \bigl(1-\nu
^M(A_j) \bigr)^{|{\mathcal L_n}|} \\
&&\qquad \leq N(\varepsilon)\max_{t\in[0,1]}t(1-t)^{|{\mathcal L_n}|}
\\
&&\qquad \leq\frac{C}{|{\mathcal L_n}| \varepsilon^{(|{M}|-1)/2}}.
\end{eqnarray*}
Combining this inequality and equality (\ref{eq:bla}), we obtain
\[
\mathbb{E}\bigl[1-S\bigl(\mathbf{X}^\star,\mathbf{X}_{(1)}^\star\bigr) | M, \mathcal L_n \bigr]
\leq\int_0^1\min\biggl(1,\frac{C}{|{\mathcal L_n}| \varepsilon
^{(|{M}|-1)/2}} \biggr)\,\dd\varepsilon.
\]
Since $|{M}|\geq4$, an easy calculation shows that there exists
$C>0$ such that
\[
\mathbb{E}\bigl[1-S\bigl(\mathbf{X}^{\star},\mathbf{X}_{(1)}^\star\bigr) | M,\mathcal L_n
\bigr]\leq
\frac{C}{|{\mathcal L_n}|^{2/(|{M}|-1)}},
\]
which leads to the desired result.
\end{pf}

\section*{Acknowledgments}
The authors are greatly indebted to Albert Ben\-ve\-nis\-te for
pointing out this problem. They also thank Kevin Bleakley and Toby
Hocking for their careful reading of the paper, and two referees and
the Associate Editor for valuable comments and insightful suggestions.

\printaddresses

\end{document}